\documentclass[10pt,psamsfonts]{amsart} 
\usepackage{amsmath}
\usepackage{amsthm}
\usepackage{amssymb}
\usepackage{amscd}
\usepackage{amsfonts}
\usepackage{amsbsy}
\usepackage{mathrsfs} 
\usepackage{setspace}
\usepackage{graphicx}
\usepackage{color}
\usepackage{calc}
\usepackage{subfigure}
\usepackage{caption}
\usepackage[font=small,labelfont=bf]{caption}
\usepackage{epsfig,afterpage}
\usepackage[dvips]{psfrag}
\usepackage[numbers,sort&compress]{natbib}
\usepackage{array}
\usepackage{enumitem}
\usepackage{multirow}

\usepackage{pgfplots}

\usepgfplotslibrary{external}
\advance\textwidth by +1.5in \advance\textheight by +1.4in
\advance\topmargin by -0.6in \advance\oddsidemargin by -0.8in
\advance\evensidemargin by -0.8in

\newcommand{\R}{\ensuremath{\mathbb{R}}}
\newcommand{\C}{\ensuremath{\mathbb{C}}}

\newcommand{\N}{\ensuremath{\mathbb{N}}}

\newtheorem {theorem} {Theorem}[section]
\newtheorem*{theoremA}{Theorem A}
\newtheorem*{theoremB}{Theorem B}

\newtheorem {lemma} [theorem]{Lemma}

\begin{document}

\title[On the Non-Existence of Isochronous Centers]
{On the Non-Existence of Isochronous Centers  in Planar Discontinuous Differential Systems}

\author[Ali Bakhshalizadeh,  Changjian Liu,  Alex C. Rezende]
{Ali Bakhshalizadeh$^1$, Changjian Liu$^2$, Alex C. Rezende$^3$}

\address{$^{1, 3}$ Departamento de Matem\'atica, Universidade Federal de S\~{a}o Carlos, S\~{a}o Carlos, Brazil.} 
\email{alibb@ufscar.br}
\email{alexcr@ufscar.br}

\address{$^{2}$ School of Mathematics (Zhuhai), Sun Yat-sen University. Zhuhai 519082, China.} 
\email{liuchangj@mail.sysu.edu.cn}



\maketitle

\begin{abstract}
The determination of whether a center is isochronous or not is a well-known problem in the qualitative theory of planar systems. In this paper, we explore planar
piecewise discontinuous differential systems characterized by a straight switching line $y=0$ and $x=0$, respectively. Our investigation reveals that such 
systems do not possess an isochronous center at the origin.
\end{abstract}

\section{Introduction and main results}\label{sec1}
One of the central challenges in the study of planar differential systems is whether a center is {\it isochronous} or not. To clarify, a center is considered 
isochronous when its {\it period function}, denoted as $T: \mathcal{S} \rightarrow \mathbb{R}$, remains constant. This period function is defined on a 
Poincar\'{e} section $\mathcal{S}$ that intersects a period annulus. This function represents the 
period of a trajectory starting from a point within $\mathcal{S}$. It is worth recalling that, for any center within a planar differential system, 
the {\it period annulus} is defined as the largest neighborhood surrounding the center that is entirely filled with periodic orbits.

In the context of smooth differential systems, the question of isochronicity in planar differential systems has received great attention from various researchers. 
For instance, please see \cite{AFG,CS,CMV,Ga,JV} and the references therein.

Regarding the piecewise discontinuous differential systems, Jules Andrade in \cite{An} proposed the following equation as a model for classical mechanical 
devices in clockmaking, which are subjected to energy loss due to friction
\begin{equation*}
\ddot{x} + \epsilon |x| \operatorname{sgn}(\dot{x}) + x = 0,
\end{equation*}
where the term $\epsilon |x| \operatorname{sgn}(\dot{x})$ models the so-called \textit{dry friction} with a damping coefficient that varies with position.
Motivated by this, the authors in \cite{MT} considered a general planar differential system of the form
\begin{equation}\label{Sys1}
(\dot{x},\dot{y})=\begin{cases}
    (y, -x-g(x)),& \text{if } y\geq 0, \\
    (y, -x+g(x)),& \text{if } y < 0,
\end{cases}
\end{equation}
where $g(x)$ is an analytic function in some neighborhood of the origin. The authors proved that if $g(x)$ is a nonlinear analytic function, then 
system \eqref{Sys1} has no isochronous center at the origin. Moreover, they studied the isochronicity of the following differential systems described by
\begin{equation}\label{Sys2}
(\dot{x},\dot{y})=\begin{cases}
    (y, -(V^{+})^{\prime}(x)),& \text{if } y\geq 0, \\
    (y, -(V^{-})^{\prime}(x)),& \text{if } y < 0,
\end{cases}
\end{equation}
and their analysis revealed that these systems exhibit isochronous behavior when $V^{\pm}(x) = a^{\pm} x^{2} + O(x^4)$, where both $a^+$ and $a^-$ 
are positive constants. Importantly, they showed that the choice of these two analytic functions can be made arbitrarily.
The authors  in \cite{LW} considered the differential systems \eqref{Sys2} where $V^{\pm}(x)$ are polynomials of even degree of the form
\begin{equation*}
\begin{aligned}
V^{+}(x)=a^{+}_{1}x^{2}+a^{+}_{2}x^{3}+\dots +a^{+}_{2n+1}x^{2n+2},\\
V^{-}(x)=a^{-}_{1}x^{2}+a^{-}_{2}x^{3}+\dots +a^{-}_{2m+1}x^{2m+2},
\end{aligned}
\end{equation*}
with $a^{+}_{i},\dots,a^{+}_{2n+1}, a^{-}_{j},\dots,a^{-}_{2m+1}$, being positive. Then, using the expansion of the period function near infinity, they
proved that the piecewise potential system \eqref{Sys2} has an isochronous center at the origin if and only if the uncoupled systems are linear.

In this paper, we begin by examining the following piecewise discontinuous differential system
\begin{equation}\label{Sys3-1}
	(\dot{x},\dot{y})=\begin{cases}
		\left( y, -(V^{+})^{\prime}(x)\right) ,& \text{if } y\geq 0, \\
		\left( (V^{-})^{\prime}(y),-x\right),& \text{if } y <0,
	\end{cases}
\end{equation}
with the switching line
\begin{equation*}
	\Sigma = \left\lbrace (x, y) \in \R^2 \,|\, y = 0\right\rbrace,
\end{equation*}
where the functions $V^{+}(x)$ and $V^{-}(y)$ are analytic within specific open regions in the plane. We impose specific conditions on the functions
\begin{equation*}
	V^{\pm} (0)=0,~(V^+)'(0) = 0,~ (V^{+})''(0) > 0, ~\text{and}~ (V^{-})'(0) < 0. 
\end{equation*}
Under the above assumptions, the upper uncoupled system exhibits a non-degenerate center at the origin, while the lower uncoupled system has an invisible quadratic tangency point at the origin.With the presence
of a non-smooth center at the origin in the coupled system, we can now present the first theorem of thepaper as follows:
\begin{theoremA}
Given the piecewise discontinuous differential system \eqref{Sys3-1}, it follows that the center at the origin does not exhibit 
isochronous behavior.
\end{theoremA}

We also consider the planar piecewise discontinuous potential system
\begin{equation}\label{Sys3}
(\dot{x},\dot{y})=\begin{cases}
   \left( y, -(V^{+})^{\prime}(x)\right),& \text{if } x\geq 0, \\
   \left( y, -(V^{-})^{\prime}(x)\right),& \text{if } x <0,
\end{cases}
\end{equation}
with the switching line
\begin{equation*}
\Sigma = \left\lbrace (x, y) \in \R^2 \,|\, x = 0\right\rbrace.
\end{equation*}
The corresponding Hamiltonian function of system \eqref{Sys3} is defined as
\begin{equation*}
H(x,y):=\frac{1}{2}y^{2}+V(x):=\begin{cases}
H^{+}(x,y)=\frac{1}{2}y^{2}+V^{+}(x), & \text{if } x \geq 0,\\
H^{-}(x,y)=\frac{1}{2}y^{2}+V^{-}(x), & \text{if } x <0,
\end{cases}
\end{equation*}
where the functions $H^{\pm}(x,y)$ are analytic within specific open regions in the plane, and $ V^{\pm} (0)=0$. Furthermore, we consider the origin as a 
non-smooth center formed by non-smooth ovals $\gamma_h\subset \lbrace H(x,y)=h\rbrace$.
Then the set of ovals $\gamma_h$ located within what we shall refer to as the period annulus, 
can be parameterized by the energy levels denoted by $h\in (0,h_{0})$, where $h_{0}\in (0,\infty)$. 

The piecewise potential system described by system \eqref{Sys3} exhibits a non-smooth center at the origin under the following conditions:
\begin{equation*}
    \begin{array}{lll}
        (i) & (V^\pm)'(0)=0, (V^\pm)''(0)>0; &  \\
        (ii) & (V^-)'(0)=\dots=(V^-)^{(2r-1)}(0)=0, (V^-)^{(2r)}(0)>0,&  \\
        & (V^+)'(0)=\dots=(V^+)^{(2s-1)}(0)=0, (V^+)^{(2s)}(0)>0;&  (r>1, s\geq 1)~\text{or}~(r\geq 1, s>1) \\
        (iii) & (V^-)'(0)=\dots=(V^-)^{(2r)}(0)=0, (V^-)^{(2r+1)}(0)<0, &  \\
        & (V^+)'(0)=\dots=(V^+)^{(2s)}(0)=0, (V^{+})^{(2s+1)}(0)>0;&(r, s\geq 1)    \\
        (iv) & (V^-)'(0)<0, (V^+)'(0)>0;&  \\
        (v) & (V^-)'(0)=0, (V^-)''(0)>0,  (V^+)'(0)>0; &\text{or} \\
        & (V^-)'(0)<0, (V^+)'(0)=0, (V^+)''(0)>0; &  \\
        (vi) & (V^-)'(0)=\dots=(V^-)^{(2r-1)}(0)=0, (V^-)^{(2r)}(0)>0, (V^+)'(0)>0; &\text{or} \\
         & (V^-)'(0)<0, (V^+)'(0)=\dots=(V^+)^{(2r-1)}(0)=0, (V^{+})^{(2r)}(0)>0; & (r> 1) \\
        (vii) & (V^-)'(0)=\dots=(V^-)^{(2r)}(0)=0, (V^-)^{(2r+1)}(0)<0, &  \\
        & (V^+)'(0)=\dots=(V^+)^{(2s-1)}(0)=0, (V^{+})^{(2s)}(0)>0; & \text{or} \\
        & (V^-)'(0)=\dots=(V^-)^{(2r-1)}(0)=0, (V^-)^{(2r)}(0)>0, &  \\
        & (V^+)'(0)=\dots=(V^+)^{(2s)}(0)=0, (V^{+})^{(2s+1)}(0)>0; & (r\geq 1, s\geq 1) \\
        (viii) & (V^-)'(0)=\dots=(V^-)^{(2r)}(0)=0, (V^-)^{(2r+1)}(0)<0, (V^+)'(0)>0; & \text{or} \\
        & (V^-)'(0)<0, (V^+)'(0)=\dots=(V^+)^{(2r)}(0)=0, (V^{+})^{(2r+1)}(0)>0; & (r\geq 1) \\
    \end{array}
\end{equation*}
where $r$ and $s$ are natural numbers. 

The authors in \cite{LW} have also conducted an investigation focused on system \eqref{Sys3} for case $(i)$ when $V^{\pm}(x)$ are the
polynomials of the form
\begin{equation*}
\begin{aligned}
V^{+}(x)=a^{+}_{1}x^{2}+a^{+}_{2}x^{3}+\dots +a^{+}_{s}x^{s+1},\\
V^{-}(x)=a^{-}_{1}x^{2}+a^{-}_{2}x^{3}+\dots +a^{-}_{l}x^{l+1},
\end{aligned}
\end{equation*}
with $ a^{+}_{i} $, $i=1,\dots,s$, and $ (-1)^{i-1}a^{-}_{i} $, $i=1,\dots,l$ being positive and $s,l\geq 1$. It is worth noting that in this particular case, 
the center of the coupled system \eqref{Sys3} is composed by the two centers from the uncoupled systems on the left and right sides.
Their findings, obtained by expanding the period function as it approaches infinity, revealed that these systems exhibit an isochronous center at the origin if and 
only if the uncoupled systems are linear. As far as we are aware, no prior research has tackled the remaining cases. In this paper, we aim to explore 
these unexplored cases.

The second main theorem of the paper is as follows:
\begin{theoremB}
Given the piecewise potential differential system \eqref{Sys3}, it follows that the centers at the origin in cases (ii) to (viii) do not exhibit 
isochronous behavior.
\end{theoremB}

The paper is organized as follows. In the second section, we begin by expanding the period function of the non-degenerate center at the origin. We also use 
the findings from paper \cite{NS}, which focus on the expansion of the period function around the origin when a differential system has a tangential center. 
With this information, we establish Theorem A. In the third section, we prove Theorem B from our paper in a step-by-step manner, addressing each case individually.

\section{Proof Of Theorem A}
This section is dedicated to the proof of the first result of this paper. 
The non-smooth center of piecewise discontinuous differential system \eqref{Sys3-1} is formed by a non-degenerate center in the upper region and a tangential 
singularity in the lower region of $y=0$. Therefore the period function $ T(x) $ of the coupled differential system \eqref{Sys3-1} is given by
\begin{equation}\label{periodfunction1}
T(x):=T^{-}(x)+T^{+}(x),
\end{equation}
where  $T^{-}(x)$ and $T^{+}(x)$ denote the periods of the lower and upper uncoupled systems, respectively.

To begin, we consider the following differential system
\begin{equation}\label{Sys4}
\begin{cases}
   \dot{x}=y, \\
    \dot{y}=-(V^{+})^{\prime}(x),
\end{cases}
\end{equation}  
where $ V^+(0)=(V^+)'(0)=0$ and $  (V^+)''(0)>0  $. Under these conditions, system \eqref{Sys4} exhibits a non-degenerate center at the origin.
In fact, the Hamiltonian function of system \eqref{Sys4} takes the following form
\begin{equation*}
H(x,y)=\frac{y^2}{2}+V^{+}(x),
\end{equation*}
where $ V^{+}(x)=\omega^{2}x^{2}+O(x^{2}) $, and all periodic orbits of System \eqref{Sys4} are contained within the level curve defined by
\begin{equation*}
\lbrace (x,y) | H(x,y)=h, h\geq 0\rbrace.
\end{equation*}
Now, consider a new variable 
\begin{equation*}
z=x\left(1+\dfrac{O(x^{2})}{\omega^{2}x^{2}}\right)^{1/2}, 
\end{equation*}
which implies $ V^{+}(z)= \omega^{2}z^{2}$. Then, we can represent $ x $ as a function of $ z $ as follows 
\begin{equation} \label{var}
x=z+\sum_{i=2}^{\infty}b_{i}z^{i}. 
\end{equation}
Hence, the period $\widetilde{T}(h)$ of the orbit  $ \gamma_h\subseteq \lbrace y^{2}/{2}+\omega^{2}z^{2}=h\rbrace $, where $y>0$, can be expressed as
\begin{equation*}
\begin{aligned}
\widetilde{T}(h)=&\int_{\gamma_h \mid_{y > 0}}\frac{1}{y}dx=\int_{\gamma_h \mid_{y > 0}}\dfrac{1+\sum_{i=1}^{\infty}(i+1)b_{i+1}z^i}{y}dz\\
=&\int_{-\frac{\sqrt{h}}{\omega}}^{\frac{\sqrt{h}}{\omega}}\dfrac{1+\sum_{i=1}^{\infty}(i+1)b_{i+1}z^i}{\sqrt{2(h-\omega^{2}z^{2})}}dz\\
=&\frac{\pi }{\sqrt{2} w}\left( 1+\sum_{i=1}^{\infty}\widetilde{T}_{2i}h^i\right),
\end{aligned}
\end{equation*}
with 
\begin{equation*}
\widetilde{T}_{2i}=\frac{\left(1+(-1)^{2 i}\right) \Gamma \left(i+\frac{3}{2}\right)}{\sqrt{\pi }\,w^{2 i}\, \Gamma (i+1)}b_{2 i+1},\qquad i=1,2,\dots
\end{equation*}
where $\Gamma$ is the gamma function. 

Noting that the function $T^{+}(x)$ can be expressed as $\widetilde{T}(h) = \widetilde{T}(V^{+}(x))$. For values of $x$ in the vicinity of zero, 
its expansion is given by
\begin{equation*}
T^{+}(x)=\frac{\pi }{\sqrt{2} w}\left(1 + \sum_{i=1}^{\infty}\widehat{T}_{2i}x^{2i}\right),
\end{equation*}
where the coefficients $\widehat{T}_{2i}$ are defined as $\widehat{T}_{2i} = \omega^{2i}\widetilde{T}_{2i}$, representing the period constants.

Now, we consider the piecewise differential system as follows
\begin{equation}\label{Sys5}
(\dot{x},\dot{y})=\begin{cases}
    \left( -(V^{-})^{\prime}(-y),-x\right),& \text{if } y\geq 0, \\
    \left( (V^{-})^{\prime}(y),-x\right),& \text{if } y <0,
\end{cases}
\end{equation}
with the following conditions $V^{-}(0)=0$ and $(V^{-})^{\prime}(0)<0$. Under these conditions, system \eqref{Sys5} exhibits a fold-fold center 
at the origin. The authors in \cite{NS} focused on determining the expansion of the period function around the origin when a differential system has 
tangential center. Specifically, they derived the power series representation of $\overline{T}(x)$ for values of $x$ near zero, which can be expressed as
\begin{equation*}
\overline{T}(x) = \sum_{i=0}^{\infty} \overline{T}_i x^i,
\end{equation*}
where
\begin{equation*}
\overline{T}_{0}=0, \quad \overline{T}_{1}=-\frac{4}{(V^-)'(0)}>0,
\end{equation*}
and the formulae for computing all the period constants are given in \cite{NS}. It is important to note that  $T^{-}(x)=\overline{T}(x)/2$.
Furthermore, using \eqref{periodfunction1} the expansion of the period function for system \eqref{Sys3-1} can be described as
\begin{equation*}
T(x)=\frac{\pi }{\sqrt{2} w}+ \sum_{i=1}^{\infty} T_i x^i,
\end{equation*}
where $T_{1}=-\frac{2}{(V^-)'(0)}>0$, $T_{2i}=\frac{\overline{T}_{2i}}{2} +\frac{\pi }{\sqrt{2} w}\widehat{T}_{2i}$, 
and $T_{2i+1}=\frac{\overline{T}_{2i+1}}{2}$ 
for $i\geq 1$. Hence, this concludes the proof for Theorem A.
\section{Proof Of Theorem B}
In this section, we will demonstrate the second theorem of this paper step by step, addressing each case individually.

In cases $(ii)$, $(iii)$, $(vi)$, $(vii)$ and $(viii)$, without loss of generality, we assume that the flatness of the left potential function $V^{-}(x)$ is equal or more 
than of the right potential function $V^{+}(x)$. For cases $(ii)$ and $(vi)$, the non-smooth center at the origin is formed by a degenerate center on the left 
side, and a center, and a tangential point on the right side, respectively. Additionally, for cases $(iii)$, $(vii)$, and $(viii)$, the 
non-smooth center at the origin is formed by a cusp on the left side, and a cusp, a center, and a 
tangential point on the right side, respectively

To describe this behavior, we define a piecewise differential system as follows
\begin{equation}\label{Sys6}
\left( \begin{array}{ll}
\dot x\\
\dot y \end{array}\right)= \begin{cases}
\left( \begin{array}{ll}
y\\
-(V^{+})^{\prime}(x) \end{array}\right), & \text{if } x \geq 0,\\
\\
\left( \begin{array}{ll}
y\\
a_{i}\, x^{2r+i-1}+O(x^{2r+i-1}) \end{array}\right), & \text{if } x <0,
\end{cases}
\end{equation}
where $i$ takes values $0$ and $1$ for cases that the left uncoupled system has a degenerate center and a cusp, respectively. 
Additionally, $a_{0}<0$ and $a_{1}>0$. The corresponding Hamiltonian function for system \eqref{Sys6} is given by
\begin{equation*}
H(x,y):=\begin{cases}
H^{+}(x,y)=\frac{1}{2}y^{2}+V^{+}(x), & \text{if } x \geq 0,\\
H^{-}(x,y)=\frac{1}{2}y^{2}-\frac{a_{i}}{2 r +i}x^{2 r +i}+O(x^{2 r +i}), & \text{if } x <0,
\end{cases}
\end{equation*}
where the origin is a non-smooth center formed by non-smooth ovals $\gamma_h:=\gamma^{-}_h\cup \gamma^{+}_h$ with
\begin{equation*}
\begin{aligned}
\gamma^{-}_h=&\left\lbrace \frac{1}{2}y^{2}-\frac{a_{i}}{2 r +i}x^{2 r +i}+O(x^{2 r +i})=h \,|\, x<0\right\rbrace,\quad
\gamma^{+}_h=&\left\lbrace \frac{1}{2}y^{2}+V^{+}(x)=h \,|\, x\geq 0\right\rbrace.
\end{aligned}
\end{equation*}

The period function for these cases can be expressed as 
\begin{equation}
T_{i}(h):=T_{i}^{-}(h)+T^{+}(h),\qquad i=0,1,
\end{equation}
where $T^{+}(h)$ represents the periods of the right uncoupled system, while $T_{i}^{-}(h)$, where $i$ takes values 0, and 1, denotes the periods of 
the left uncoupled system for cases that it has a degenerate center and a cusp, respectively. The left uncoupled system exhibits different behaviors for cases with 
a degenerate center and a cusp, depending on whether we are in cases $(ii)$, $(iii)$, $(vi)$, $(vii)$, or $(viii)$. 
In cases $(ii)$ and $(iii)$, the left uncoupled system has a degenerate center at the origin, while in cases $(vi)$, $(vii)$, and $(viii)$, it has a cusp at the origin.

Assuming $h=\rho^{2r+i}$, for $i=0, 1$, we can express $y$ in the left uncoupled system as follows
\begin{equation*}
y=\pm \sqrt{2\left( \rho^{2r+i}+\frac{a_{i}}{2 r +i}x^{2r+i}+O(x^{2 r +i})\right) },
\end{equation*}
subsequently, we can compute $ T_{i}^{-}(\rho) $ as
\begin{equation*}
\begin{aligned}
T_{i}^{-}(\rho)=&2\int_{\sqrt[2r+i]{-\frac{2r+i}{a_{i}}}\rho+O(\rho)}^{0}\dfrac{dx}{\sqrt{2\left( \rho^{2r+i}+\frac{a_{i}}{2 r +i}x^{2r+i}+O(x^{2 r +i})\right) }}\\
=&\sqrt{\frac{2}{\rho^{2(r-1)+i}}}
\int_{\sqrt[2r+i]{-\frac{2r+i}{a_{i}}}+O(\rho^{0})}^{0}\dfrac{ds}{\sqrt{1+\frac{a_{i}}{2 r +i}s^{2r+i}+O(\rho^{0}s^{2 r +i})}},
\end{aligned}
\end{equation*}
where, at the last step, we applied the change of variable $x = \rho s$, yielding
\begin{equation*}
\lim_{\rho\rightarrow 0}T_{i}^{-}(\rho)=\infty.
\end{equation*}

Consequently, the period function $T_{i}(h)$ for $i=0, 1$ tends towards infinity near the origin. Therefore, in cases $(ii)$, $(iii)$, $(vi)$, $(vii)$ and $(viii)$, 
the non-smooth center 
at the origin does not exhibit isochronous behavior.

In case $(iv)$, both $H^+$ and $H^-$ tangentially intersect the switching line $\Sigma$ at the point $(0,0)$. This point $(0,0)$ serves as an invisible 
quadratic tangency for both vector fields, resembling an invisible fold-fold point. Consequently, based on Theorem $A$ in \cite{NS}, there is no existence of an 
isochronous center in this particular case. 

In case $(v)$, without loss of generality, we assume that system \eqref{Sys3} exhibits a non-degenerate center on the left side and an invisible quadratic tangency on the right side. Our goal is to 
determine the expansion of the period function around $x=0$. We define the period function of the coupled system \eqref{Sys3} in this case as follows
\begin{equation}\label{periodfunction}
T(y):=T^{-}(y)+T^{+}(y),
\end{equation}
where $T^{-}(y)$ and $T^{+}(y)$ denote the periods of the left and right uncoupled systems, respectively.

First, we consider the following differential system
\begin{equation}\label{Sys7}
\begin{cases}
   \dot{x}=y, \\
    \dot{y}=-(V^{-})^{\prime}(x),
\end{cases}
\end{equation}  
where the potential function is defined as $V^{-}(x) = \omega^2 x^2 + O(x^2)$, which implies the presence of a center at the origin.
By applying the variable transformation $(x, y, t) \mapsto (x,\sqrt{2}\,\omega y, \frac{1}{\sqrt{2}\omega}\, t)$, we can simplify the system \eqref{Sys7} to the 
following form
 \begin{equation}\label{Sys8}
\begin{cases}
   \dot{x}=y:=f(x,y), \\
    \dot{y}=-\dfrac{1}{2\omega^{2}}(V^{-})^{\prime}(x):=g(x,y),
\end{cases}
\end{equation} 
 where it possesses an equilibrium point at the origin with the linear part having eigenvalues $\pm i$.
Then, system \eqref{Sys8} in the complex coordinate $ z=x+iy $ can be written as 
\begin{equation}\label{Sys9}
\dot{z}=h(z,\bar{z}),\quad \dot{\bar{z}}=\bar{h}(z,\bar{z}),
\end{equation}
where 
\begin{equation*}
h(z,\bar{z})=f\left( \frac{z+\bar{z}}{2},\frac{z-\bar{z}}{2i}\right) +i g\left( \frac{z+\bar{z}}{2},\frac{z-\bar{z}}{2i}\right),
\end{equation*}
is analytic function in some neighbourhood $ \nu^{*} $ of the origin in $\C^2$. Note that it is straightforward to show that 
$h(z,\bar{z})=i z+O(|z,\bar{z}|^{2})$.
Given that the second component of equation \eqref{Sys9} is essentially the complex conjugate of the first component, our primary focus is directed towards 
the following equation
\begin{equation}\label{Sys10}
\dot{z}=h(z,\bar{z}).
\end{equation}

By applying near the identity changes of variables, as the spirit of normal form transformations, system \eqref{Sys10} in complex coordinate can be simplified to 
\begin{equation*}
\dot{z}=i z+\sum_{j=1}^{N}\left( \alpha_{2j+1}z(z\bar{z})^{j}+i \beta_{2j+1}z(z\bar{z})^{j}\right)+O(2N+3),
\end{equation*}
for any $N\in \N$, where $\alpha_{2j+1},\beta_{2j+1}\in \R$, please see \cite{AFG}. In the polar coordinate the normal form can be expressed 
as 
\begin{equation}\label{polar}
\begin{cases}
\dot{r}=r^{2N+3}R(r,\theta),\\
\dot{\theta}=1+\beta_{3}r^{2}+\beta_{5}r^{4}+\dots+\beta_{2N+1}r^{2N}+r^{2N+2}\Theta (r,\theta),
\end{cases}
\end{equation}
 where the functions $ R(r,\theta) $ and $ \Theta (r,\theta) $
are analytical in $r$ and $2\pi$-periodic in $\theta$. Consequently, we obtain the following equality
\begin{equation*}
\dfrac{dr}{d\theta}=\dfrac{r^{2N+3}R(r,\theta)}{1+\beta_{3}r^{2}+\beta_{5}r^{4}+\dots+\beta_{2N+1}r^{2N}+r^{2N+2}\Theta (r,\theta)}.
\end{equation*}

In the following lemma, we will now present the period function for system \eqref{Sys7}, specifically when the condition $x < 0$ holds.
The proof of this lemma follows the same approach as that of Theorem 2.1 in \cite{AFG}, with the additional consideration that 
$\frac{\pi}{2} < \theta < \frac{3}{2}\pi$ and a subsequent reversion to the original variable.
\begin{lemma}
Consider system \eqref{Sys7} under the condition $x< 0$. Then, the period function near $y=0$ is expressed as
\begin{equation*} 
T^-(y) = \frac{\pi}{\sqrt{2}\omega}\left(1 + \sum_{i=1}^{\infty}\frac{1}{(\sqrt{2}\omega)^i}\widetilde{T}_{i}y^{i}\right), 
\end{equation*}
with
\begin{equation*}
\widetilde{T}_{2i-1} = 0, \quad \widetilde{T}_{2i} = \sum_{\substack{n_1 + \dots + n_l = 2i\\ n_j \text{ even},\ l \geq 1}} (-1)^l \beta_{n_1 + 1} 
\dots \beta_{n_l + 1},
\end{equation*}
where the coefficients $\widetilde{T}_{i}$ represent the period constants.
\end{lemma}

Now, we consider the piecewise differential system as follows
\begin{equation}\label{Sys11}
(\dot{x},\dot{y})=\begin{cases}
    (y, -(V^{+})^{\prime}(x),& \text{if } x\geq 0, \\
    (y, (V^{+})^{\prime}(-x),& \text{if } x <0,
\end{cases}
\end{equation}
subject to the conditions that $V^{+}(0) = 0$ and $(V^{+})^{\prime}(0) > 0$. Under these conditions, system \eqref{Sys11} exhibits a fold-fold center at 
the origin. Using the results from \cite{NS}, we can express the power series of the period function for values of $y$ near zero as
\begin{equation*}
\overline{T}(y) = \sum_{i=0}^{\infty} \widehat{T}_i y^i,
\end{equation*}
where
\begin{equation*}
\widehat{T}_{0}=0, \quad \widehat{T}_{1}=\frac{4}{(V^+)'(0)}>0,
\end{equation*}
and the formulae for computing all the period constants are given in \cite{NS}. Therefore, we have that $T^{+}(y)=\overline{T}(y)/2$.
Furthermore, using \eqref{periodfunction} the expansion of the period function for system \eqref{Sys3} in case $(v)$ can be described as
\begin{equation*}
T(y)=\frac{\pi}{\sqrt{2}\omega}+ \sum_{i=1}^{\infty} T_i y^i,
\end{equation*}
where $T_{1}=\frac{2}{(V^+)'(0)}>0$, $T_{2i}=\frac{\widehat{T}_{2i}}{2} +\frac{\pi}{(\sqrt{2}\omega)^{2i+1}}\widetilde{T}_{2i}$, 
and $T_{2i+1}=\frac{\widehat{T}_{2i+1}}{2}$ 
for $i\geq 1$. Therefore, this shows that there is no center at the origin exhibiting isochronous behavior.

\section*{Acknowledgments}
The first author is supported by FAPESP grant, process number 2022/07822-5. The second author......
The third author is partially support by FAPESP grant, process number 2019/09385-9.

\end{document}